\documentclass[11pt]{amsart}
\usepackage[margin=1in]{geometry}

\usepackage{amssymb}
\usepackage{amsthm}
\usepackage{amsmath}
\usepackage{mathrsfs}
\usepackage{amsbsy}
\usepackage[all]{xy}
\usepackage{bm}
\usepackage{hyperref}
\usepackage{tikz}
\usepackage{array}
\usepackage{enumerate}
\usepackage{xcolor}
\usepackage{bbm}
\usepackage{comment}
\usepackage{dynkin-diagrams}
\usepackage{ifthen}
\usepackage{thmtools, thm-restate}

\hypersetup{colorlinks=true, citecolor=lavender, linkcolor=Magenta,urlcolor=Magenta}

\definecolor{Magenta}{rgb}{1,0,0.6}
\definecolor{lavender}{rgb}{0.4,0,1}

\usepackage{hhline}
\allowdisplaybreaks
\usepackage[noadjust]{cite}
\usepackage{asymptote}

\usepackage{caption}
\usepackage{tabu}
\usepackage{diagbox}
\usepackage[noabbrev,capitalise,nameinlink]{cleveref}
\crefname{conjecture}{Conjecture}{Conjectures}

\newcommand{\C}{\mathbb{C}}
\newcommand{\R}{\mathbb{R}}
\newcommand{\E}{\mathbb{E}}

\newcommand{\diff}{\mathop{}\!\mathrm{d}}
\newcommand{\dz}{\diff z}
\newcommand{\dw}{\diff w}

\DeclareMathOperator{\des}{des}
\DeclareMathOperator{\Var}{Var}

\newtheorem{theorem}{Theorem}[section]

\newtheorem{lemma}[theorem]{Lemma}

\theoremstyle{definition}

\newtheorem{remark}[theorem]{Remark}

\begin{document}

\title[]{Central limits from generating functions}

\author[]{Mitchell Lee}
\address[]{Department of Mathematics, Harvard University, Cambridge, MA 02138, USA}
\email{mitchell@math.harvard.edu}

\begin{abstract}
    Let $(Y_n)_n$ be a sequence of $\R^d$-valued random variables. Suppose that the generating function \[f(x, z) = \sum_{n = 0}^\infty \varphi_{Y_n}(x) z^n,\] where $\varphi_{Y_n}$ is the characteristic function of $Y_n$, extends to a function on a neighborhood of $\{0\} \times \{z : |z| \leq 1\} \subset \R^d \times \C$ which is meromorphic in $z$ and has no zeroes. We prove that if $1 / f(x, z)$ is twice differentiable, then there exists a constant $\mu$ such that the distribution of $(Y_n - \mu n) / \sqrt{n}$ converges weakly to a normal distribution as $n \to \infty$.

    If $Y_n = X_1 + \cdots + X_n$, where $(X_n)_n$ are i.i.d.\@ random variables, then we recover the classical (Lindeberg--L\'evy) central limit theorem. We also prove the 2020 conjecture of Defant that if $\pi_n \in \mathfrak{S}_n$ is a uniformly random permutation, then the distribution of $(\des (s(\pi_n)) + 1 - (3 - e) n) / \sqrt{n}$ converges, as $n \to \infty$, to a normal distribution with variance $2 + 2e - e^2$.
\end{abstract}

\maketitle

\section{Introduction}\label{section:intro}

For any positive integer~$d$, let $\langle \cdot, \cdot \rangle \colon \R^d \times \R^d \to \R$ denote the standard bilinear form given by $\langle x, y \rangle = x_1 y_1 + \cdots + x_d y_d$. For any $\R^d$-valued random variable~$Y$, let $\varphi_Y \colon \R^d \to \C$ denote the corresponding \emph{characteristic function}, given by \[\varphi_Y(\omega) = \E[\exp(i \langle Y, \omega\rangle)],\] where~$\E$ denotes expected value \cite[Chapter~3.3]{MR3930614}. We denote the covariance matrix of~$Y$ by $\Var[Y]$. 

We denote partial derivatives using a subscript; for example, $g_{x_1 z}(x, z)$ denotes \[\frac{\partial}{\partial x_1} \frac{\partial}{\partial z} g(x, z).\]

For any real, symmetric, and positive semidefinite matrix $\Sigma \in \R^{d \times d}$, let $\mathcal{N}(0, \Sigma)$ denote the multivariate normal distribution with mean~$0$ and covariance matrix~$\Sigma$ \cite[Chapter~3]{MR0560319}.

The main theorem of this article is the following central limit theorem. Like the classical (Lindeberg--L\'evy) central limit theorem, it states that a particular sequence of random variables converges in distribution to a normally distributed random variable.

\begin{theorem}\label{thm:clt}
Let~$d$ be a positive integer, and let $Y_0, Y_1, Y_2, \ldots$ be a sequence of $\R^d$-valued random variables. Suppose that there exists a function $g \colon U \to \C$, where $U$ is an open neighborhood of $\{0\} \times \{z : |z| \leq 1\} \subset \R^d \times \C$, such that
\begin{enumerate}[(i)]
    \item $g(x, z)$ is holomorphic as a function of $z$ for any fixed $x$;
    \item $g$ is twice differentiable;
    \item for all $(x, z) \in U$ with $|z| < 1$, we have \begin{equation}\label{eq:expected-value-power-series}\sum_{n = 0}^\infty \varphi_{Y_n}(x) z^n = \frac{1}{g(x, z)}.\end{equation}
\end{enumerate}

For all $j, k$ with $1 \leq j, k \leq d$, define \[\mu_j = i g_{x_j}(0, 1)\] and \[\Sigma_{j, k} = g_{x_j x_k}(0, 1) - i (\mu_j g_{x_k z}(0, 1) + \mu_k g_{x_j z}(0, 1)) + \mu_j \mu_k.\] Then, we have $\lim_{n \to \infty} \E[Y_n]/n = \mu$ and $\lim_{n \to \infty} \Var[Y_n]/n = \Sigma$. Moreover, $Z_n = (Y_n - \mu n) /\sqrt{n}$ converges in distribution, as $n \to \infty$, to $Z \sim \mathcal{N}(0, \Sigma)$.
\end{theorem}

\cref{thm:clt} is similar to the following central limit theorem proved by Bender in 1973, which has been shown to be useful throughout analytic combinatorics \cite{MR1039294,MR1041444,MR1181718,MR1217750}.

\begin{theorem}[{\cite[Theorem~1]{MR375433}}] \label{thm:clt-bender-richmond}
    For all $n, k \geq 0$, let $a_n(k)$ be a nonnegative integer. Suppose that for any fixed $n$, only finitely many of the $a_n(k)$ are nonzero.

    For all $n$, define $Y_n$ to be the $\mathbb{N}$-valued random variable that takes the value $k$ with probability \[\frac{a_n(k)}{\sum_{i = 0}^\infty a_n(i)}.\] 

    Let \[f(z,w) = \sum_{n = 0}^\infty\sum_{k = 0}^\infty a_n(k) w^k z^n\] be the bivariate generating function of the sequence $a_n(k)$. Suppose that there exist a function $A(s)$ continuous and nonzero near $0$, a function $r(s)$ with bounded third derivative near $0$, a nonnegative integer $m$, and $\epsilon, \delta > 0$ such that
    \[\left(1 - \frac{z}{r(s)}\right) f(z, e^s) - \frac{A(s)}{1 - z/r(s)}\] is analytic and bounded for $|s| < \epsilon$ and $|z| < |r(0)| + \delta$. Define \[\text{$\mu = -\frac{r'(0)}{r(0)}$ \qquad and \qquad $\sigma^2 = \mu^2 - \frac{r''(0)}{r(0)}$}.\] If $\sigma \neq 0$, then $Z_n = (Y_n - \mu n)/\sqrt{n}$ converges in distribution, as $n \to \infty$, to $Z \sim \mathcal{N}(0, \sigma)$.
\end{theorem}
Bender and Richmond also proved a multivariate generalization of \cref{thm:clt-bender-richmond} in 1983 \cite[Corollary~1]{MR0700034}. Both \cref{thm:clt} and \cref{thm:clt-bender-richmond} are central limit theorems whose hypotheses refer to a particular power series in the variable $z$. The primary difference between the two theorems is that \cref{thm:clt} refers to the series \[\sum_{n = 0}^\infty \varphi_{Y_n}(x) z^n,\] in which the coefficient of $z^n$ is the characteristic function of the real-valued random variable $Y_n$. In contrast, \cref{thm:clt-bender-richmond} refers to the series \[\sum_{n = 0}^\infty\sum_{k = 0}^\infty a_n(k) w^k z^n,\] in which the coefficient of $z^n$ is $\sum_{k = 0}^\infty a_n(k) w^k$, which is a scaled form of the probability generating function of the $\mathbb{N}$-valued random variable $Y_n$. Another difference is that the hypotheses of \cref{thm:clt} can be checked more directly, which we will use in the proof of \cref{corollary:defant} below.

In \cref{section:proof}, we will prove \cref{thm:clt}.

In \cref{section:applications}, we will show how \cref{thm:clt} easily implies the Lindeberg--L\'evy central limit theorem:
\begin{restatable}[{\cite[Theorem~3.10.7]{MR3930614}}]{corollary}{llclt}\label{corollary:llclt}
    Let $(X_n)_n$ be $\R^d$-valued i.i.d.\@ random variables such that $\E\left[|X_1|^2\right] < \infty$. Let $\mu = \E[X_1]$ and $\Sigma = \operatorname{Cov}[X_1]$. Then \[\frac{X_1 + \cdots + X_n - \mu n}{\sqrt{n}}\] converges in distribution, as $n \to \infty$, to $Z \sim \mathcal{N}(0, \Sigma)$.
\end{restatable}

Then, we will prove the following 2020 conjecture of Defant.
\begin{restatable}[{\cite[Conjecture~7.11]{MR4384616}}]{corollary}{defant}\label{corollary:defant}
    Let $\des \colon \mathfrak{S}_n \to \mathbb{N}$ denote the function that counts the descents of a permutation, and let $s \colon \mathfrak{S}_n \to \mathfrak{S}_n$ denote West's stack-sorting map. If $\pi_n \in \mathfrak{S}_n$ is a uniformly random permutation, then \[\frac{\des (s(\pi_n)) + 1 - (3 - e) n}{\sqrt{n}}\] converges in distribution, as $n \to \infty$, to $Z \sim \mathcal{N}(0, 2 + 2e - e^2)$.
\end{restatable}
\begin{remark}
    Clearly, the expression $\des (s(\pi_n)) + 1$ appearing in the statement of \cref{corollary:defant} can be replaced by $\des (s(\pi_n))$, but we use $\des (s(\pi_n)) + 1$ to match the conventions of \cite{MR4384616}.
\end{remark}

\section*{Acknowledgements}
The author thanks Colin Defant for helpful correspondence and Lutz Mattner for providing valuable feedback on the article.

\section{Proof of Theorem~\ref{thm:clt}}\label{section:proof}

\begin{proof}[Proof of \cref{thm:clt}]
    Let $D = \{z : |z| \leq 1\}$ be the closed unit disc. Since~$D$ is compact, we may replace $U$ with a smaller convex set $U_1 \times U_2$, where $U_1$ is an open neighborhood of $0 \in \R^d$ and $U_2$ is an open neighborhood of $D \subseteq \C$.

    Let 
    \begin{equation}\label{eq:f}
        f(x, z) = \frac{1}{g(x, z)} = \sum_{n = 0}^\infty \varphi_{Y_n}(x) z^n.
    \end{equation} We start by proving that it is possible to take the derivative of $f(x, z)$ at $x = 0$ term by term.
    
    \begin{lemma}\label{lemma:differentiable}
        For all $n$, the characteristic function $\varphi_{Y_n}(x)$ is twice differentiable at $x = 0$. Moreover, we have \[f_{x_j}(0, z) = \sum_{n = 0}^\infty \left(\frac{\partial}{\partial x_j}\varphi_{Y_n}\right)(0) z^n\] for $1 \leq j \leq d$ and $|z| < 1$. 
    \end{lemma}
    \begin{proof}
        By the hypotheses of \cref{thm:clt}, $g(x, z)$ is twice differentiable and nonvanishing for $(x, z) \in U$ with $|z| < 1$. Therefore, $f(x, z)$ is twice differentiable for $(x, z) \in U$ with $|z| < 1$ as well.

        Fix any $z$ with $|z| < 1$ and let $r = (1 + |z|) / 2$. By the Cauchy integral formula and \eqref{eq:f}, we have \begin{equation}\label{eq:phi-y_n}\varphi_{Y_n}(x) = \frac{1}{2 \pi i} \oint_{|w| = r} \frac{f(x, w)}{w^{n + 1}}\dw \end{equation} for all $x \in U_1$. Since the integrand $f(x, w) / w^{n + 1}$ of \eqref{eq:phi-y_n} is twice differentiable as a function of $x$ at $x = 0$, the function $\varphi_{Y_n}(x)$ is twice differentiable as a function of $x$ at $x = 0$ as well.
        
        Moreover, we may differentiate both sides of \eqref{eq:phi-y_n} using differentiation under the integral sign to conclude that
        \[\left(\frac{\partial}{\partial x_j}\varphi_{Y_n}\right)(0) = \frac{1}{2 \pi i} \oint_{|w| = r} \frac{f_{x_j}(0, w)}{w^{n + 1}}\dw\] for $1 \leq j \leq d$.
        The term $f_{x_j}(0, w)$ is bounded for $|w| = r$. Therefore, by the dominated convergence theorem,
        \begin{align*}
            \sum_{n = 0}^\infty \left(\frac{\partial}{\partial x_j}\varphi_{Y_n}\right)(0) z^n &= \sum_{n = 0}^\infty \left(\frac{1}{2 \pi i} \oint_{|w| = r} \frac{f_{x_j}(0, w)}{w^{n + 1}}\dw\right) z^n \\
            &= \frac{1}{2 \pi i} \oint_{|w| = r} \left(f_{x_j}(0, w)\sum_{n = 0}^\infty \frac{z^n}{w^{n + 1}} \right)\dw \\
            &= \frac{1}{2 \pi i} \oint_{|w| = r} \frac{f_{x_j}(0, w)}{z - w} \dw \\
            &= f_{x_j}(0, z),
        \end{align*}
        where the last equality follows from the Cauchy integral formula.
    \end{proof}

    Let us now show that~$\mu_j$ is real for $1 \leq j \leq d$. By \cref{lemma:differentiable} and \cite[Section~XV.4]{MR270403}, the random variable $Y_n$ has a finite expectation for all $n$, with components given by \[i\E[(Y_n)_j] = \left(\frac{\partial}{\partial x_j}\varphi_{Y_n}\right)(0).\] We may take the derivative of $g(x, z) = 1/f(x, z)$ with respect to $x_j$ on both sides and substitute $x = 0$. This yields
    \[g_j(0, z) = -\frac{\sum_{n = 0}^\infty \left(\frac{\partial}{\partial x_j} \varphi_{Y_n}\right)(0) z^n}{\left(\sum_{n = 0}^\infty \varphi_{Y_n}(0) z^n\right)^2} = -\frac{\sum_{n = 0}^\infty i \E[(Y_n)_j] z^n}{\left(\frac{1}{1 - z}\right)^2} = -(1-z)^2\sum_{n = 0}^\infty i \E[(Y_n)_j] z^n.\]
    Therefore,
    \[i g_j(0, z) = (1 - z)^2 \sum_{n = 0}^\infty \E[(Y_n)_j] z^n.\]
    As $z \to 1^-$ through real numbers, the left-hand side of this equation approaches $\mu_j$ and the right-hand side is always real. Therefore, $\mu_j$ is real.

    It follows that $Z_n = (Y_n - \mu n) / \sqrt{n}$ is an $\R^d$-valued random variable for all~$n$. We will now prove the last statement of the theorem that $Z_n$ converges in distribution to $Z \sim \mathcal{N}(0, \Sigma)$. By L\'evy's convergence theorem \cite[Theorem~3.3.17]{MR3930614}, it suffices to show that for all $\omega \in \R^d$, we have
    \begin{equation}\label{eq:desired-limit}\lim_{n \to \infty}\varphi_{Z_n}(\omega) = \exp\left(-\frac{1}{2}\langle \omega, \Sigma \omega \rangle\right).\end{equation}

    Broadly, we will prove \eqref{eq:desired-limit} by writing $\varphi_{Z_n}(\omega)$ in terms of $\varphi_{Y_n}\left(\omega / \sqrt{n}\right)$. To estimate the latter, we will split the series \eqref{eq:expected-value-power-series} into a principal part and an analytic part, and show that only the principal part contributes meaningfully to the value of $\varphi_{Y_n}(\omega / \sqrt{n})$.

    Observe that by substituting $x = 0$ into \eqref{eq:expected-value-power-series}, we obtain $g(0, z) = 1 - z$. Therefore,~$g$ only has one zero on the compact set $\{0\} \times D$. It is at $(x, z) = (0, 1)$, with $g(0, 1) = 0$ and $g_z(0, 1) = -1 \neq 0$. Therefore, by the implicit function theorem, there exist a neighborhood $V_1 \subseteq U_1$ of $0 \in \R^d$, a neighborhood $V_2 \subseteq U_2$ of $D \subseteq \C$, and a twice differentiable function $b \colon V_1 \to \C$ such that for all $(x, z) \in V_1 \times V_2$, we have
    \begin{equation*}\label{eq:implicit-function}
        \text{$g(x, z) = 0$\quad if and only if\quad $z = b(x)$.}
    \end{equation*}
    Since $b(0) = 1$, we may also assume, by replacing $V_1$ with a smaller open set, that~$b$ does not vanish on $V_1$.

    For any fixed $x \in V_1$, the function $f(x, z) = 1 / g(x, z)$ is meromorphic on $V_2$ and its only pole is at $z = b(x)$. Hence, we may remove its pole by subtracting the principal part. Explicitly, the function
    \[h(x, z) = f(x, z) - \frac{a(x)}{1 - z / b(x)} = \sum_{n = 0}^\infty \left(\varphi_{Y_n}(x) - \frac{a(x)}{(b(x))^n}\right) z^n\]
    extends analytically from $V_2 \setminus \{r(x)\}$ to $V_2$, where
    \[a(x) = -\frac{g_z(x, b(x))}{b(x)}.\]

    Now, we proceed in a manner similar to Bender and Richmond \cite[Corollary~1]{MR0700034}. Since $V_2$ is an open neighborhood of $D$, it contains the closed ball $\{z : |z| \leq r\}$ for some $r > 1$. For all~$n$ and~$x$, the coefficient of $z^n$ in the series $h(x, z)$, which we denote $[z^n]h(x, z)$, can be computed using the Cauchy integral formula:
    \[[z^n]h(x, z) = \frac{1}{2 \pi i}\oint_{|z| = r} \frac{h(x, z)}{z^{n + 1}} \dz.\]
    Therefore,
    \[|[z^n]h(x, z)| \leq r^{-n} \sup_{|z| = r} |h(x, z)| = O(r^{-n}),\]
    where the constant hidden by the~$O$ notation is uniform for~$x$ in any compact set.

    It follows that \begin{equation}\label{eq:yn}\varphi_{Y_n}(x) = [z^n] f(x, z) = \frac{a(x)}{(b(x))^n} + [z^n]h(x, z) = \frac{a(x)}{(b(x))^n} + O(r^{-n}),\end{equation} where, again, the constant hidden by the $O$ is uniform for $x$ in any compact set.

    Let us now turn to \eqref{eq:desired-limit}. Fix $\omega \in \R^d$. For all~$n$, we have
    \begin{align}
        \varphi_{Z_n}(\omega) &= \E\left[\exp\left(i \left\langle \frac{Y_n - \mu n}{\sqrt{n}}, \omega\right\rangle\right)\right] \nonumber \\
        &= \frac{\E\left[i \left\langle Y_n, \frac{\omega}{\sqrt{n}}\right\rangle \right]}{\exp(i \langle\mu, \omega\rangle \sqrt{n})} \nonumber \\
        &= \frac{a\left(\frac{\omega}{\sqrt{n}}\right)}{\left(b\left(\frac{\omega}{\sqrt{n}}\right)\right)^n\exp(i \langle\mu, \omega\rangle \sqrt{n})} + O(r^{-n})\label{eq:zn-omega}
    \end{align}
    where we used \eqref{eq:yn} in the third equality.

    Now, we analyze the denominator of \eqref{eq:zn-omega}. Since~$b$ is twice differentiable and $b(0) = 1$, the function $\log(b(x))$ has a second order Taylor expansion for $x \to 0$ \cite[Theorem~12.14]{MR0344384}. We may compute the coefficients by differentiating the equation $g(x, b(x)) = 0$ using the chain rule. This yields \[\log (b(x)) = -i\langle \mu, x \rangle + \frac{1}{2} \langle x, \Sigma x\rangle + o(|x|^2).\]
    Therefore, recalling that $\omega$ is fixed,
    \begin{align*}
        &\left(b\left(\frac{\omega}{\sqrt{n}}\right)\right)^n\exp(i \langle\mu, \omega\rangle \sqrt{n})\\
        &= \exp\left(n \log \left(b\left(\frac{\omega}{\sqrt{n}}\right)\right) + i \langle\mu, \omega\rangle \sqrt{n}\right) \\
        &= \exp\left(n\left(-i \left\langle \mu, \frac{\omega}{\sqrt{n}}\right\rangle + \frac{1}{2} \left\langle \frac{\omega}{\sqrt{n}}, \Sigma\frac{\omega}{\sqrt{n}}\right\rangle + o\left(\frac{1}{n}\right) \right) + i \langle\mu, \omega\rangle \sqrt{n}\right) \\
        &= \exp\left(\frac{1}{2}\langle \omega, \Sigma \omega\rangle + o(1)\right).
    \end{align*}
    Substituting into \eqref{eq:zn-omega}, we have
    \begin{align*}
        \lim_{n \to \infty}\varphi_{Z_n}(\omega) &= \frac{\displaystyle\lim_{n \to \infty} a\left(\frac{\omega}{\sqrt{n}}\right)} {\displaystyle\lim_{n \to \infty} \left(b\left(\frac{\omega}{\sqrt{n}}\right)\right)^n\exp(i \langle\mu, \omega\rangle \sqrt{n}) } \\
        &= \frac{a(0)}{\exp\left(\frac{1}{2} \langle\omega, \Sigma \omega\rangle\right)} \\
        &= \exp\left(-\frac{1}{2} \langle\omega, \Sigma \omega\rangle\right),
    \end{align*}
    proving \eqref{eq:desired-limit}. This completes the proof of the last statement of the theorem.

    It remains to prove that $\lim_{n \to \infty} \E[Y_n]/n = \mu$ and $\lim_{n \to \infty} \Var[Y_n]/n = \Sigma$, which can be done using the following computations:
    \[
        \lim_{n \to \infty} \frac{\E[Y_n]}{n} = \lim_{n \to \infty} \frac{\E[\sqrt{n}Z_n + \mu n]}{n} = \left(\lim_{n \to \infty}n^{-1/2}\right) \left(\lim_{n \to \infty} \E[Z_n]\right) + \mu = \mu 
    \]
    and
    \[
        \lim_{n \to \infty} \frac{\Var[Y_n]}{n} = \lim_{n \to \infty} \frac{\Var[\sqrt{n}Z_n + \mu n]}{n} = \lim_{n \to \infty} \Var[Z_n] = \Var[Z] = \Sigma. \qedhere
    \]
\end{proof}

\section{Applications}\label{section:applications}
We now prove \cref{corollary:llclt,corollary:defant} to demonstrate the utility of \cref{thm:clt}.
\llclt*
\begin{proof}
    Let $Y_n = X_1 + \cdots + X_n$. We have
    \[\sum_{n = 0}^\infty \varphi_{Y_n}(x) z^n = \sum_{n = 0}^\infty (\varphi_{X_1}(x))^n z^n = \frac{1}{1 - \varphi_{X_1}(x) z}.\]
    Now, apply \cref{thm:clt} with $g(x, z) = 1 - \varphi_{X_1}(x) z$. This is clearly holomorphic in $z$, and it is twice differentiable because $\E\left[|X_1|^2\right] < \infty$ \cite[Section~XV.4]{MR270403}.
\end{proof}
\defant*
\begin{proof}
    Following Defant, let \begin{equation}\label{eq:F}F(y, z) = \frac{y}{2} (-1 - y z + \sqrt{1 - 4z + 2yz + y^2z^2})\end{equation} where we choose the branch of the square root that evaluates to $1$ as $z \to 0$.
    
    Let \begin{equation}\label{eq:F-hat}\hat{F}(y, z) = \sum_{m, n = 0}^\infty F_{m, n} \frac{y^m z^n}{n!},\end{equation} where $F_{m, n}$ is the coefficient of $y^mz^n$ in $F(y, z)$. (The power series $\hat{F}$ can alternatively be defined by \[\hat{F}(y, z) = \mathcal{L}^{-1}\{F(y, 1/t)/t\}(z),\] where $\mathcal{L}^{-1}$ is the inverse Laplace transform with respect to the variable $t$.)
    
    We have \cite[Theorem~7.8]{MR4384616}
    \begin{equation}\label{eq:des}
        \sum_{n = 1}^\infty \left(\sum_{\pi \in \mathfrak{S}_{n-1}} y^{\des(s(\pi)) + 1}\right)\frac{z^n}{n!} = - \log(1 + \hat{F}(y, z)).
    \end{equation}
    Differentiating with respect to $z$ and substituting $y = e^{ix}$, we find that the generating function of the characteristic functions $\varphi_{\des(s(\pi_n)) + 1}$ is the following:
    \[
        \sum_{n = 0}^\infty \varphi_{\des(s(\pi_n)) + 1}(x) z^n = -\frac{\hat{F}_z(e^{ix}, z)}{1 + \hat{F}(e^{ix}, z)}.
    \]
    Let
    \begin{equation}\label{eq:defant-g}
        g(x, z) = -\frac{1 + \hat{F}(e^{ix}, z)}{\hat{F}_z(e^{ix}, z)}
    \end{equation}
    be the reciprocal of this generating function. We now check that $g(x, z)$ has an analytic continuation that satisfies the conditions of \cref{thm:clt}.

    First, observe that by \eqref{eq:F}, $F(y, z) / y$ can be written as a power series in $z$ and $yz$. Therefore, the coefficient $F_{m, n}$ is zero if $m > n + 1$. It is also easy to check using Darboux's lemma \cite[Theorem~1]{MR1003845} that $F_{m, n}$ is bounded by an exponential function of $n$. Hence, the series \eqref{eq:F-hat} converges everywhere, so $\hat{F}$ is entire (on $\C^2$). Therefore, the numerator and denominator of \eqref{eq:defant-g} are analytic (real-analytic in the first variable and complex-analytic in the second variable) as a function of $(x, z) \in \R \times \C$.

    One may easily compute $F(1, z) = -z$ by substituting $y = 1$ in \eqref{eq:F}. It follows that $\hat{F}(1, z) = -z$ as well, so $\hat{F}_{z}(1, z) = -1$. Therefore, there is a neighborhood $U$ of $\{0\} \times \{z : |z| \leq 1\}$ on which $\hat{F}_z(e^{ix}, z)$ does not vanish. Then $g$ is analytic on $U$.

    By \cref{thm:clt}, there exist constants $\mu, \sigma$ such that 
    \[\frac{\des(s(\pi_n)) + 1 - \mu n}{\sqrt{n}}\]
    converges in distribution, as $n \to \infty$, to $Z \sim \mathcal{N}(0, \sigma)$. We can compute $\mu = 3 - e$ and $\sigma = 2 + 2e - e^2$, either by using the definition of $\mu$ and $\sigma$ in the statement of \cref{thm:clt}, or by using the formulas for $\E[\des(s(\pi_n)) + 1]$ and $\operatorname{Var}[\des(s(\pi_n)) + 1]$ \cite[Theorems~7.9~and~7.10]{MR4384616}.
\end{proof}

\bibliographystyle{plain}
\bibliography{main}

\end{document}